\documentclass[12pt]{paper}

\usepackage{amscd}
\usepackage[titletoc,toc,title]{appendix}

\usepackage{amsthm}
\usepackage{thmtools}

\usepackage[headings]{fullpage}
\usepackage{mathrsfs}
\usepackage{amsmath,amsfonts,amssymb,amsthm}
\usepackage[all]{xy}
\usepackage{amscd}

\usepackage[dvips]{graphics,color}
\usepackage{hyperref}

\newtheorem*{thma}{Theorem 1}
\newtheorem*{thmb}{Theorem 2}

\def\noi{\noindent}

\def\pf{\noi{\bf Proof.\ \,}}

\def\eop{{$\square$}}

\def\QQ{{\mathbb Q}}

\def\ZZ{{\mathbb Z}}

\def\mp{{\mathfrak p}}
\def\mq{{\mathfrak q}}

\def\Spec{{\mathrm{Spec}}}

\def\Spc{{\mathrm{Spc}}}

\def\perf{{\mathrm{perf}}}

\def\al{ {\cal A}  }

\def\th{\theta}

\begin{document}

\newtheorem{thm}{Theorem}[section]
\newtheorem{prop}[thm]{Proposition}
\newtheorem{lem}[thm]{Lemma}
\newtheorem{coro}[thm]{Corollary}

\theoremstyle{definition}
\newtheorem{de}[thm]{Definition}
\newtheorem{nota}[thm]{Notation}

\newtheorem{rem}[thm]{Remark}
\newtheorem{ex}[thm]{Example}

\begin{center}

{\bf \ Classifying space of subcategories and its application}
\medskip

Yong Liu \footnote{Email address is \emph{yongliue@gmail.com}}

Department of Mathematics and Statistics

Auburn University, Alabama, U.S.

\end{center}

\begin{abstract}

For a collection of subcategories satisfying a fixed set of conditions, for example thick subcategories of a triangulated category, we define a topological space called classifying space of subcategories. We show that this space classifies various prime subcategories in the sense that they bijectively correspond to the closed subsets of the classifying space. Many well-known results of subcategory classification fit into this framework. An example which cannot be classified by a topological space in the above sense is also given.

\end{abstract}

\tableofcontents

\section{Introduction}

The classification of subcategories of various types is a fundamental problem with different background. For examples, Gabriel's paper~\cite{Gabriel} showed that in the category $R$-$\mathrm{mod}$ of $R$-modules the Serre subcategories are classified via the specialization closed subsets of the prime spectrum $\Spec~R$. The $t$-structures~\cite{BBD} corresponding to the subcategories, aisles, are also popular and many papers attempt to give an invariant such as~\cite{Stan} and~\cite{Alonso2}. There are many other famous results claiming that the target subcategories are one to one correspondent to the closed (or dually open) subsets of certain spectrum, such as Neeman's classification of (co)localizing subcategories in~\cite{Neeman} and~\cite{Neeman2}. Balmer gave a classification of the radical thick tensor ideals in a tensor triangulated cateogry in~\cite{Balmer}. The classification of subcategories can also be used for reconstruction, especially for schemes. This is achieved by Rosenberg~\cite{Rosenberg} and Rouquier~\cite{Rouquier}, by using appropriate spectra.

In our paper, we study this phenomenon in a more general context, in which a collection $\Phi$ of subcategories we are interested in forming a complete lattice is considered. Then a space $K(\Phi)$ is defined, consisting of points $P_C$ such that $C\in \Phi$ cannot be generated as proper subcategories of the same type (see Definition~\ref{points} and Definition~\ref{type}). This space $K(\Phi)$ however is not in general a topological space in a natural way unless we require the lattice $\Phi$ also be distributive. This problem is resolved by restricting our attention into a subset of (generally) prime elements (that is, whenever $C\leq \bigvee_{i\in I}C_i$, $C\leq C_i$ for some $i$, where $I$ is finite or arbitrary accordingly), denoted by $K_{p}(\Phi)$ (or $K_{gp}(\Phi)$).  Therefore, we obtain our main result as Theorem~\ref{closedprime}, after a brief discussion on connection among those classifying spaces.

\begin{thma}\label{introthm}
Assume that $\Phi$ is a complete lattice of subcategories of certain type in a category $\al$. Then there is an isomorphism of lattices
\[
\{\text{closed subsets of}~K_{gp}(\Phi)\}\stackrel{\sim}{\rightarrow}\{\text{g-primely generated subcategories of certain type in}~\al\}.
\]
\end{thma}

Interestingly, thanks to lattice duality, there are many similar notions to (generally) prime objects, and if we ignore the condition of point representation, a family of parallel results to Theorem~\ref{closedprime} such as Theorem~\ref{classsification} and  Corollary~\ref{cortoclassification} are obtained, which can be considered as various criterions of subcategories that are classified by a topological space.

\begin{thmb}
Let $\Phi$ be a complete lattice of subcategories of certain type in a category $\al$. Then the subcategories generated by join primes are classified by the closed subsets of $K^{\vee}_p(\Phi)$. The dual statements and their complete versions for the meet primes and irreducibles also hold accordingly.
\end{thmb}

This provides a framework for many famous results including Kanda's classification of Serre subcategories in~\cite{Kanda}, Neeman's classification~\cite{Neeman} of thick subcategories in $D_{\perf}(D)$ and also the classification of localizing subcategories in an noetherian stable homotopy category~\cite{HoveyPS}, by showing that each subcategory is g-primely generated, respectively. And we show for the last example as a demonstration in Section~\ref{cmpclassical}.

We point out that also as another comparison, in P. Johnstone's book~\cite{JohnstoneP} on Stone spaces, he defined a similar notion using lattices, which was used to develop the theory of classifying subcategories in a derived category by J. Kock and W. Pitsch recently in~\cite{Kock}.

As other classifications of various subcategories establish a correspondence from either closed or open subsets of a topological space to the target subcategories, we mention in Section~\ref{exceptional} an exceptional example of nullity class in an abelian category, which cannot be classified by a topological space (see Definition~\ref{classifies}) but corresponds to the closed subsets with one extra condition, namely extension-closed, see Theorem 6.8 in~\cite{liu2}.



\section{Subcategories and lattices}

There is a natural way to equip a collection of subcategories with a structure of complete (distributive) lattice. For details in lattice theory, see~\cite{DaveyPriestley} and~\cite{JohnstoneP} for example. We introduce in this section the concept that subcategories are classified by a topological space in general.

\begin{de}

A lattice $(\Phi,\leq)$ is a partially ordered set with two commutative binary operations defined via the partial order $\leq$, called \textit{join} $\vee$ and \textit{meet} $\wedge$ (or \textit{supremum} and \textit{infimum}, resp.), such that any finte subset of $\Phi$ has both join and meet. In particular, the join of the empty set denoted by 0 is called the \textit{bottom element}, and the meet of the empty set denoted by 1 is called the \textit{top element}.


We call $\Phi$ a \textit{distributive lattice} if additionally one of the following properties holds for any elements $a,b,c$ in $\Phi$: (1) $a\vee(b\wedge c)=(a\vee b)\wedge (a\vee c)$; (2) $a\wedge(b\vee c)=(a\wedge b)\vee (a\wedge c)$. In fact, that one holds implies the other, see Lemma 1.5 in~\cite{JohnstoneP}. A lattice $\Phi$ is \textit{complete} if it has either arbitrary joins or arbitrary meets for any subset of $\Phi$.


\end{de}

There are two typical examples of lattices coming from either a category or a topological space.

\begin{ex}\label{2example}
(1) Assume $\al$ is an abelian category such that the collection $\Phi_{\mathrm{tor}}$ of its torsion classes forms a set. Equipped with the inclusion as partial order, the intersection as meet, and generating a torsion class as join, then $\Phi_{\mathrm{tor}}$ becomes a complete lattice. However, it is not always distributive as we will see in Section~\ref{exceptional}.

(2) Let $X$ be a topological space and $\Phi_X$ the set of closed subsets of $X$. Then $\Phi_X$ is a complete distributive lattice with the inclusion as partial order, the intersection as meet and the closure of their union as join.
\end{ex}


\begin{de}
Let $\Phi_1,\Phi_2$ be partially ordered sets. A \textit{homomorphism} of partially ordered sets is a map $f:\Phi_1\rightarrow \Phi_2$ that preserves the partial orders. Such an $f$ becomes \textit{isomorphism} if it has an inverse homomorphism. For lattices, we require additionally that the map respects both joins and meets, and we define isomorphism of lattices in a similar manner.
\end{de}

It is easy to show that any isomorphism of partially ordered sets is automatically an isomorphism of lattices once the partially ordered sets are also lattices.




The following notion in Definition~\ref{type} is introduced for our convenience, so that we use it to represent any lattice of subcategories in a given category. For instance, in an abelian category the lattice of Serre subcategories, and in a triangulated category the lattice of localizing subcategories and so on.

\begin{de}\label{type}
Let $\al$ be a category. The collection $\Phi_A$ of \textit{subcategories of type} $A$ (or roughly, \textit{subcategories of certain type}) in $\al$ is the collection of subcategories in which every subcategory satisfies a fixed (finite) set of conditions.
\end{de}

The next definition in fact describes a common phenomenon in classification of subcategories.

\begin{de}\label{classifies}
Let $\al$ be a category. The collection $\Phi$ of subcategories of certain type in $\al$ is \textit{classified by a topological space} $X$ if there is an isomorphism of partially ordered sets $f:\Phi_X\stackrel{\sim}{\rightarrow}\Phi$, where $\Phi_X$ is the lattice of closed subsets of $X$.
\end{de}


For example, Gabriel's result~\cite{Gabriel} in 1962 says that the collection of Serre subcategories in $R$-$\mathrm{mod}$ is classified by the prime spectrum $\Spec R$ of the ring $R$ with specialization closed topology. Neeman~\cite{Neeman} showed in 1992 that the collection of localizing subcategories in $D(R)$ is classified by $\Spec R$ with discrete topology. Another one given by Balmer in 2005 says that the collection of radical thick tensor ideals in a tensor triangulated category $\al$ is classified by the Hochster's dual $(\Spc \al)^h$ of Balmer's spectrum, consisting of the prime thick tensor ideals of $\al$, see~\cite{Balmer} and~\cite{Hochster}. Also, Benson et al. proved~\cite{Benson} in 2011 that the collection of localizing tensor ideals in the homotopy category $K(\mathrm{Inj} kG)$ is classified by the spectrum $\Spec H^\ast(G,k)$ of the cohomology of a finite group $G$ with discrete topology.

We will connect some classical results mentioned above to our frame work in Section~\ref{cmpclassical} once the notion of prime subcategory in a classifying space of subcategories is defined. For convenience, the set theoretical issue is ignored and we consider lattices instead of general collections of subcategories.



\section{Classifying space of subcategories}

In this section, we construct for any lattice $\Phi$ a topological space $K(\Phi)$. The lattice $\Phi$ comes from either a category or a topological space as in Example~\ref{2example}. We use the capital letters such as $C$ to denote an element in $\Phi$.

\begin{de}\label{points}
Let $\Phi$ be a lattice. Assign to each element $C\in\Phi$ a pairing $(C, C^o)$, called \textit{a point} and denoted by $P_C$, if $C\neq C^o$, where $C^o$ is the supremum of all elements strictly below, namely
\[
C^o=\bigvee\{C'~|~C'\lneqq C,C'\in\Phi\}.
\]
Notice that the collection of points in $\Phi$ may be empty. The collection $K(\Phi)=\{P_C~|~C\in \Phi\}$ is called \textit{the classifying space associated to} $\Phi$.

\end{de}




According to the lattices, we have the following properties about the points $P_C$, roughly claiming that every point has a representative.

\begin{lem} \label{closure}
(1) Let $\Phi_A$ be the lattice of subcategories of type A in a category. Then for each $P_C\in K(\Phi_A)$, there is an object $x\in C$ such that $C=\langle x\rangle$ is the intersection of all $C\in\Phi_A$ containing $x$.

(2) Let $\Phi_X$ be the lattice of closed subsets of a topological space $X$. Then for each $P_C\in K(\Phi_X)$, there is an element $x\in C$ such that $C=\overline{\{x\}}$, where $\overline{\{x\}}$ denotes the closure of $\{x\}$ in $X$.
\end{lem}


\pf We show (2), and (1) follows similarly. Suppose $C-\bigvee_{\substack{C'\lneqq C\\C'\in\Phi}}C'\neq \emptyset$ and $x$ is an element in this difference. If $\overline{\{x\}}\neq C$, then the closure has to be one of the $C'$s in $\Phi$ properly contained in $C$. Therefore, $x\in\bigvee_{\substack{C'\lneqq C\\C'\in\Phi}}C'$, a contradiction.~\eop


For every $C\in\Phi$, denote by $K(C)=\{P_{C'}\in K(\Phi)~|~C'\leq C\}$ a subset of $K(\Phi)$. Clearly this collection contains the whole space $K(\Phi)$ by taking the top element. However, it does not contain the empty set, and we add it artificially. We thus expect the candidate $\{K(C)\}_{C\in\Phi}$ to give a topology on $K(\Phi)$.

Recall that a topological space is $T_0$ (or \textit{Kolmogorov}) if for any two distinct points there is an open neighborhood containing one point but not the other. The \textit{$T_0$-quotient} $KQ(X)=X/\sim$ of a space $X$ is defined as a quotient space by the equivalence relation that $x\sim y$ if and only if $\overline{\{x\}}=\overline{\{y\}}$.

\begin{thm}\label{kquotient}
Suppose $\Phi$ is a complete distributive lattice. Then the collection $\{K(C)\}_{C\in\Phi}$ defines a topology of closed sets on $K(\Phi)$, making $K(\Phi)$ a $T_0$-space.

Moreover, suppose $\Phi$ is the lattice $\Phi_X$ of closed subsets of a topological space $X$. If the closure of any $x\in X$ represents a point in $K(\Phi_X)$, then $K(\Phi_X)$ is homeomorphic to the $T_0$-quotient of $X$.
\end{thm}

\pf By definition, $\bigcap_{i\in I} K(C_i)=K(\bigwedge_{i\in I} C_i)$ for an arbitrary index $I$ and $C_i\in \Phi$. Now suppose $C_1,C_2\in\Phi$, it suffices to show $K(C_1\vee C_2)\leq K(C_1)\bigcup K(C_2)$. That is, for any $C\leq C_1\vee C_2$ representing a point $P_C\in K(\Phi)$, either $C\leq C_1$ or $C\leq C_2$. Suppose not, then the distributivity implies that there is a nontrivial decomposition $C=C\wedge(C_1\vee C_2)=(C\wedge C_1)\vee(C\wedge C_2)$, contradicting the fact that $P_C$ is a point. Furthermore, let $P_{C_1}, P_{C_2}\in K(\Phi)$ be distinct points. The complement $K(C_1)^c$ would give an open neighborhood of $P_{C_2}$ but not containing $P_{C_1}$. Otherwise, $P_{C_2}\in K(C_1)$ implies $C_2\lneqq C_1$, which allows us to choose $K(C_2)^c$ as a separation instead. Hence $K(\Phi)$ is $T_0$.


Suppose $\Phi_X$ is the lattice of closed sets of a space $X$. There is a natural map
\[
\phi:KQ(X)\rightarrow K(\Phi_X)
\]
by assigning an equivalence class $[x]$ to $P_{\overline{\{x\}}}$, which is well-defined by our assumption. Clear that $\phi$ is injective.  Now for any $P_C\in K(\Phi_X)$, Lemma~\ref{closure} implies $\phi([x])=P_{\overline{\{x\}}}=P_C$ for some $x\in C-\bigvee_{\substack{C'\lneqq C\\C'\in\Phi}}C'$. Hence $\phi$ is surjective. Finally, it is straightforward to check that such $\phi$ is closed and continuous by observing that $\phi[C]=K(C)$ and $p^{-1}\phi^{-1}K(C)=C$, where $[C]$ denotes a closed subset of $KQ(X)$ and $p:X\rightarrow KQ(X)$ is the projection.
~~~\eop

It is not difficult to see that for a space $X$, the classifying space $K(\Phi_X)$ can also be defined via the lattice of open subsets of $X$ accordingly, and they turn out to be homeomorphic.


\begin{rem}\label{primeness}
In $K(\Phi_X)$ for any space $X$, every point $P_C$ is represented by a single point $x\in X$ by Lemma~\ref{closure}. It follows that $\{K(C)\}_{C\in\Phi_X}$ automatically gives a topology of closed subsets on $K(\Phi_X)$. In particular, each point $P_C$ is prime as we will see later. However, it is not true in general for the lattices of subcategories of certain type.

\end{rem}




We show that the construction of classifying space $K(\Phi)$ respects isomorphisms.

\begin{lem}\label{func2}
Let $f:\Phi_1\rightarrow \Phi_2$ be an isomorphism of lattices. Then

\noindent(1) $C\in\Phi_1$ represents a point in $K(\Phi_1)$ if and only if $f(C)$ represents a point in $K(\Phi_2)$;

\noindent(2) $f$ induces a homeomorphism $K(\Phi_1)\stackrel{\approx}{\rightarrow} K(\Phi_2)$.
\end{lem}

\pf Statement (2) follows from (1), and (1) holds because $C^o\neq C$ if and only if $f(C)^o\neq f(C)$.~~~\eop

\begin{itemize}

\item \textit{Examples of classifying space of subcategories }

\end{itemize}

We end this section by giving three examples of $K(\Phi)$ for lattice $\Phi$ of subcategories of certain type. First, let $\Phi_r$ denote the lattice of \textit{replete} subcategories of a category $\al$, that is, each subcategory is closed under isomorphisms of objects. We assume the ambient category is essentially small for convenience, so that a lattice of subcategories makes sense.

\begin{prop}\label{repletecat}

The classifying space $K(\Phi_r)$ is homeomorphic to the space of isomorphism classes of objects with discrete topology.

\end{prop}

\pf Take an element $C\in\Phi_r$ and objects $x,y\in C$. We claim that either $x\cong y$ or $C=C^o$. In fact, suppose $x\ncong y$. Then there is a nontrivial decomposition $C=\langle x\rangle_r \cup (\langle x\rangle_r )^c$, where $\langle x\rangle_r $ is the smallest replete subcategory containing $x$ with its complement $(\langle x\rangle_r )^c$ in $C$, and notice that that $(\langle x\rangle_r )^c$ is also a replete subcategory. Conversely, each isomorphism class of an objet $x\in\al$ gives a replete subcategory $\langle x\rangle_r$ which contains no proper replete subcategories in $\langle x\rangle_r$. Hence it represents a point $P_{\langle x\rangle_r}$ in $K(\Phi_r)$. Therefore, the underlying set of $K(\Phi_r)$ corresponds to the set of isomorphism classes of objects. Since both $\langle x\rangle_r$ and its complement $\langle x\rangle_r^c$ in $\al$ are replete, and also $K(\langle x\rangle_r)^c= K(\langle x\rangle_r^c)$, every singleton $K(\langle x\rangle_r)=\{P_{\langle x\rangle_r}\}$ is open. Thus this topology is discrete.~~~\eop


Furthermore, consider in a \textit{Krull-Schmidt category} $\al$ (that is, every object can be written as a coproduct of finite many indecomposables), the lattice $\Phi_{\coprod}$ of replete subcategories that are also closed under retracts and finite coproducts.

\begin{prop}
The classifying space $K(\Phi_{\coprod})$ has points corresponding to the isomorphism classes of indecomposable objects.
\end{prop}

\pf By Lemma~\ref{closure}, choose a representative $x\in\al$ for a point $P_C$. If there is a decomposition say $x=x_1\coprod x_2$, then either $\langle x_1\rangle_{\coprod}=\langle x\rangle_{\coprod}$ or $\langle x_2\rangle_{\coprod}=\langle x\rangle_{\coprod}$, where $\langle a\rangle_{\coprod}$ represents the smallest subcategory in $\Phi_{\coprod}$ containing $a$. In fact, if it is not the case, then $\langle x\rangle_{\coprod}$ cannot represent a point in $K(\Phi_{\coprod})$. Therefore, each representative $x$ of $P_C$ can be chosen to be indecomposable. The conclusion then follows from the fact that $\langle x\rangle_{\coprod}$ contains no proper subcategories of the same type, for any indecomposable $x$. ~~~\eop

Lastly, consider the lattice $\Phi_s$ of Serre subcategories in an abelian category $\al$. An object $x$ is \textit{monoform} if it contains no subquotient of $x$ as a subobject of $x$, other than the trivial cases. Denote by $\langle x\rangle_s$ the smallest Serre subcategory containing $x$. See~\cite{Kanda} for detail.

\begin{prop}
Suppose $\al$ is a noetherian abelian category. Then the points in $K(\Phi_s)$ correspond to a subset of the isomorphism classes of monoform objects in $\al$.
\end{prop}

\pf First we show that the representative $x$ of every point $P_C=P_{\langle x\rangle_s}$ can be chosen to be monoform. In fact, suppose $x$ is not monoform, then by Proposition 2.2 in~\cite{Kanda} it has a monoform subobject $x_1$ such that either $\langle x\rangle_s=\langle x_1\rangle_s$ or $\langle x\rangle_s=\langle x/x_1\rangle_s$. Since otherwise $\langle x\rangle_s=\langle x_1\rangle_s\vee\langle x/x_1\rangle_s$ implies that $\langle x\rangle_s$ does not represent a point in $K(\Phi_s)$. If $\langle x\rangle_s=\langle x_1\rangle_s$ or $\langle x\rangle_s=\langle x/x_1\rangle_s$ such that $x/x_1$ is monoform, then we are done. If $\langle x\rangle_s=\langle x/x_1\rangle_s$ but $x/x_1$ is not monoform, then we can replace $x$ by $x/x_1$ and repeat the procedure again, say a subobject $x_2$ of $x$ with a monomorphism $x_1\hookrightarrow x_2$. Thus we can obtain a chain of subobjects of $x$
\[
x_1\hookrightarrow x_2\hookrightarrow x_3\hookrightarrow\cdots\hookrightarrow x,
\]
which has to stop after finitely many stages, since $\al$ is noetherian. Then $x_n=x_{n+1}=\cdots$ and $\langle x\rangle_s=\langle x/x_1\rangle_s=\cdots=\langle x/x_n\rangle_s$. Since the last term $x/x_n$ contains no proper subobject, it has to be monoform.~~~\eop

The converse statement is not necessarily true because a monoform object may not represent a point. For example, in the abelian category $\al$ of finite dimensional representations of type $A_2:\stackrel{1}{\circ}\longleftarrow\stackrel{2}{\circ}$ over any field, there are three indecomposables which are also monoform, namely $P_1, P_2$ and $S_2$ with a short exact sequence $0\rightarrow P_1\rightarrow P_2\rightarrow S_2\rightarrow 0$. However, this gives a nontrivial decomposition of the Serre subcategory $\al=\langle P_2\rangle_s=\langle P_1\rangle_s\vee\langle S_2\rangle_s$.

\section{Generally prime objects and classification}\label{gprime}

By introducing the concept of (generally) prime elements of a lattice $\Phi$ of subcategories, we obtain a classification of these subcategories, requiring no distributivity of $\Phi$. We also compare various classifying spaces of subcategories in this section.


\begin{de}
Let $\al$ be a category and $\Phi$ a lattice of subcategories in $\al$ of certain type. A subcategory $C\in\Phi$ or a point $P_C\in K(\Phi)$ is called \textit{prime} (\textit{generally prime} or simply \textit{g-prime}, resp.) if whenever $C\leq C_1\vee C_2$ with $C_i\in\Phi$ we have either $C\leq C_1$ or $C\leq C_2$ ($C\leq\bigvee_{i\in I} C_i$ implies $C\leq C_i$ for some $i\in I$ with $I$ an arbitrary index set, resp.).
\end{de}

In particular, if $P_C$ is represented by a single object $x\in\al$, that the category $C$ is prime is equivalent to that whenever $x\in C_1\vee C_2$ we have either $x\in C_1$ or $x\in C_2$. We denote by
\[
K_p(\Phi)=\{P_C\in K(\Phi)~|~C~\text{is prime}\}
\]
the set of prime points.

Now as before choose a collection of subsets $\{K_p(C)\}_{C\in\Phi}$ in which
\[
K_p(C)=\{P_{C'}\in K_p(\Phi)~|~C'\leq C\}.
\]

\begin{prop}
Let $\al$ be a category and $\Phi$ a lattice of subcategories in $\al$ of certain type. Then the collection $\{K_p(C)\}_{C\in\Phi}$ forms a topology of closed subsets of $K_p(\Phi)$, which is $T_0$.
\end{prop}

\pf Since $\emptyset,\al\in\Phi$, we have $K_p(\emptyset)=\emptyset$ and $K_p(\al)=K_p(\Phi)$. It is clear that $\bigcap K_p(C_i)=K_p(\bigwedge C_i)$ since one can take the prime points of both sides from the equality $\bigcap K(C_i)=K(\bigwedge C_i)$. We need to show that $K_p(C_1)\vee K_p(C_2)=K_p(C_1\vee C_2)$ for each $C_i\in\Phi$. In fact, take the primes of both sides from $K(C_1)\vee K(C_2)\leq K(C_1\vee C_2)$ to obtain one containment. For the converse, take $P_C\in K_p(C_1\vee C_2)$. Then $C\leq C_1\vee C_2$ so that $C\leq C_1$ or $C\leq C_2$ by primeness of $C$. The same argument of Proposition~\ref{kquotient} applies to the second assertion.~~~\eop

Similar to Lemma~\ref{func2} we also have

\begin{lem}
Let $f:\Phi_1\rightarrow \Phi_2$ be an isomorphism of lattices. Then it induces a homeomorphism $f_p:K_p(\Phi_1)\stackrel{\approx}{\rightarrow} K_p(\Phi_2)$.
\end{lem}

We define $K_{gp}(\Phi)$ similarly
\[
K_{gp}(\Phi)=\{P_C\in K(\Phi)~|~C~\text{is generally prime}\}
\]
as the set of generally prime points and show without difficulty that it has a topology of closed subsets given by
\[
K_{gp}(C)=\{P_{C'}\in K_{gp}(\Phi)~|~C'\leq C\}.
\]

\begin{prop}
The space $K_{gp}(\Phi)$ is a $T_0$-space with the topology of closed subsets $\{K_{gp}(C)\}_{C\in\Phi}$. Moreover, any isomorphism of lattices $f:\Phi_1\stackrel{\sim}{\rightarrow}\Phi_2$ induces a homeomorphism $f_{gp}:K_{gp}(\Phi_1)\stackrel{\approx}{\rightarrow} K_{gp}(\Phi_2)$.
\end{prop}

\begin{ex}
The whole category is usually prime but not g-prime. Consider the derived category $\al=D^b_{fg}(\QQ)$ of bounded chain complexes over the rationals $\QQ$ with finitely generated homologies, and $\Phi$ the lattice of nullity classes. Since the generators are of the form $\Sigma^i\QQ$ for some $i\in \ZZ$ and every nullity class becomes $\langle\Sigma^i\QQ\rangle_n\cap\al$ with $\langle\Sigma^i\QQ\rangle_n$ the nullity class generated in $D(\QQ)$, the category $\al$ as a nullity class is prime because the nullity classes are truncated from below. However, it is not g-prime since we have
\[
\al=\bigvee_{i\in\ZZ}\langle\Sigma^i\QQ\rangle_n\cap\al
\]
with every $\langle\Sigma^i\QQ\rangle_n\cap\al$ properly contained in $\al$. See~\cite{Stan} for detail.

\end{ex}

Next we give some observations on their connection among various classifying spaces for the lattice of .

\begin{prop}
Let $\Phi$ be a lattice from either a topological space $X$ or subcategories of certain type.

(1) If $X$ has specialization-closed topology, then $K_{gp}(\Phi_X)=K(\Phi_X)$;

(2) For any well-ordered index set $I$, if union $\bigcup_{i\in I}C_i$ of any chain
$\cdots\leq C_{i-1}\leq C_i\leq C_{i+1}\leq\cdots$ in $\Phi$ remains in $\Phi$, then $K_{gp}(\Phi)=K_p(\Phi)$.
\end{prop}

\pf (1) Suppose an arbitrary union of closed subsets of $X$ is closed, and $P_C\in K(\Phi_X)$. By Lemma~\ref{closure}, $C=\overline{\{x\}}$ for some $x\in X$. Thus if $C=\bigvee_{i\in I}C_i=\bigcup_{i\in I}C_i$ with $C_i$ closed, we have $x\in C_i$ so that $C\leq C_i$ for some $i$.


(2) Assume $C\in\Phi$ is prime and $C\leq\bigvee_{i\in I} C_i$ with $C_i\in \Phi$. Then we form a chain $C'_i:=\bigvee_{j\in I, j\leq i}C_i$, so that $\bigvee_{i\in I} C_i=\bigvee_{i\in I} C'_i=\bigcup_{i\in I} C_i'$ by assumption. Suppose $C$ has a representative $x$ and let $i\in I$ be the least element such that $x\in C'_i$. Then we claim $x\in C_i$. Indeed, if $i$ has a successor, then $C_i'=C_{i-1}'\vee C_i$ so that the primeness of $C$ implies that $x\in C_i$ since $x\notin C'_{i-1}$. For those $i$ which has no successor, we have
\[
C_i'=(\bigvee_{j\in I,j<i}C_j)\bigvee C_i=(\bigcup_{j\in I, j<i}C_j')\bigvee C_i
\]
in which $\bigcup_{j\in I,j<i}C_j'\in\Phi$ again by assumption. Hence the primeness of $C$ implies either $x\in\bigcup_{j\in I, j<i}C_j'$ or $x\in C_i$. However, that $x\in\bigcup_{j\in I, j<i}C_j'$ would imply $x\in C_j'$ for some $j<i$, contradicting the minimality of $i$. Therefore, $x\in C_i$ as well.~~~\eop

Let $\Phi$ be a complete lattice of certain subcategories of a category $\al$. Notice that every generally prime subcategory necessarily represents a point in $K(\Phi)$ by definition, since otherwise $C=C^o=\bigvee_{C'\lneqq C,C'\in\Phi} C'$ implies $C\leq C'\lneqq C$ for some $C'\in\Phi$ by g-primeness, which is absurd. We call a subcategory of $\al$ is \textit{generated by (g-)primes} or \textit{(g-)primely generated} if it can be written as a join of subset of (g-)prime objects.


\begin{thm}\label{closedprime}
Assume that $\Phi$ is a complete lattice of subcategories of certain type in $\al$. Then there is a bijection
\[
\th:\{K_{gp}(C)\}_{C\in\Phi}\stackrel{\sim}{\rightarrow}\{C\in\Phi~|~C=\widehat{C}\}:\xi
\]
where $\widehat{C}=\bigvee_{\substack{C'\leq C\\ C'~\text{g-prime}}}C'$. In particular, there is an isomorphism of lattices
\[
\{\text{closed subsets of}~K_{gp}(\Phi)\}\stackrel{\sim}{\rightarrow}\{\text{g-primely generated subcategories of certain type}\}.
\]
\end{thm}


\pf Define $\th$ as
\[
\th(K_{gp}(C))=\bigvee_{P_{C'}\in K_{gp}(C)}C':=D.
\]
Notice that every g-prime $C\in\Phi$ represents a point in $K(\Phi)$. Therefore, $\th(K_{gp}(C))=\widehat{C}=D$. Then by comparing the index sets, the fact $C'\leq \widehat{C}=D$ implies
\[
D=\bigvee_{P_{C'}\in K_{gp}(C)}C'=\bigvee_{\substack{C'\leq C\\ C'~\text{g-prime}}}C'\leq\bigvee_{\substack{D'\leq D\\ D'~\text{g-prime}}}D'=\widehat{D}\leq D.
\]
Hence $\th$ is well-defined because $\widehat{D}=D$. For the inverse map, set $\xi(C)=K_{gp}(C)$. The order preserving property follows immediately from the definitions of $\th$ and $\xi$.

Next, for every $C$ with $C=\widehat{C}$,
\[
\th\xi(C)=\th(K_{gp}(C))=\widehat{C}=C.
\]
Also for every $K_{gp}(C)$ with $C\in\Phi$,
\[
\xi\th(K_{gp}(C))=\xi(\widehat{C})=K_{gp}(\widehat{C})=K_{gp}(C),
\]
where the last equality holds because $\widehat{C}\leq C$ implies $K_{gp}(\widehat{C})\subseteq K_{gp}(C)$, while that every point in $K_{gp}(C)$ is represented by some g-prime $E\in\Phi$ contained in $C$ implies $E\leq \bigvee\limits_{\substack{C'\leq C\\ \text{g-prime}~C'\in\Phi}}C'=\widehat{C}$ so that $K_{gp}(\widehat{C})\supseteq K_{gp}(C)$.~~~\eop

Theorem~\ref{closedprime} can be restated as the following, in terms of Definition~\ref{classifies}.

\begin{coro}
For a complete lattice $\Phi$ of subcategories of a certain type in a category $\al$, these subcategories are classified by $K_{gp}(\Phi)$ if they are g-primely generated.
\end{coro}

\section{Irreducible vs. prime, and functoriality}\label{irred}

In this section, we briefly discuss in a lattice the irreducible and the prime elements in various versions. In fact, the two notions coincide if the lattice is distributive shown in Lemma~\ref{irredandprime}. Then we show that our construction $K_{gp}(\Phi)$ of a complete lattice $\Phi$ is functorial in a non-obvious way, shedding light on the point-free construction of Johnstone~\cite{JohnstoneP} presented in Section II.1.3, where the point in $\Phi$ is described as a homomorphism $p:\Phi\rightarrow\{0,1\}$ of lattices with $\{0,1\}$ the lattice consisting of the bottom and the top elements.

By \textit{dualization} of a lattice $(\Phi,\leq,\vee,\wedge)$, we refer to the lattice $(\Phi,\geq,\wedge,\vee)$ with the same underlying set but with reversed order and join, meet interchanged. In other words, we turn the graph of $\Phi$ upside down.

\begin{de}
Let $\Phi$ be a (complete) lattice and  $C, C_i\in\Phi$ for $i\in I$. Then

(1) $C$ is \textit{join irreducible} if $C=C_1\vee C_2$ $\Longrightarrow$ $C=C_1$ or $C=C_2$;

($1'$) $C$ is \textit{completely join irreducible} if $C=\bigvee C_i$ $\Longrightarrow$ $C=C_i$ for some $i$;

(2) $C$ is \textit{join prime} if $C\leq C_1\vee C_2$ $\Longrightarrow$ $C\leq C_1$ or $C\leq C_2$;

($2'$) $C$ is \textit{completely join prime} if $C\leq\bigvee C_i$ $\Longrightarrow$ $C\leq C_i$ for some $i$.

Their dual notions in $\Phi$, \textit{(completely) meet irreducible} and \textit{(completely) meet prime} elements, are defined respectively.

\end{de}

The notion of completely join prime element is identified with that of generally prime element in Section~\ref{gprime}. Also notice that completely join (meet resp.) irreducible implies finite join (meet resp.) irreducible, and similarly for primeness.

\begin{lem}\label{irredandprime}
Let $\Phi$ be a lattice and $C\in\Phi$. Then

(1) if $C$ is join prime, then it is join  irreducible;

(2) if additionally $\Phi$ is distributive, then that $C$ is join irreducible implies it is join prime.

The statements of complete lattices with the dual notions involved are also true.
\end{lem}

\pf We only show for the finite versions and the complete cases are similar.

(1) Call $C\in\Phi$ \textit{iprime} if $C=C_1\vee C_2$ implies $C\leq C_1$ or $C\leq C_2$. Clearly, prime implies iprime. Now suppose $C$ is iprime. Then $C=C_1\vee C_2$ implies that $C\leq C_1$ or $C\leq C_2$, so that $C_1\leq C\leq C_1$ or $C_2\leq C\leq C_2$. Hence iprime implies irreducible.

(2) Suppose $C\leq C_1\vee C_2$. Then $C=(C\wedge C_1)\vee(C\wedge C_2)$ by distributivity, so that either $C=C\wedge C_1$ or $C=C\wedge C_2$ by irreducibility of $C$. Hence $C\leq C_1$ or $C\leq C_2$.~~~\eop

The problem of functoriality becomes clear if we use the point-free description, see Section II.1.3 in Johnstone~\cite{JohnstoneP}. Consider
\[
K_p^{\wedge}(\Phi)= \{C\in\Phi~|~C~\text{is meet prime}\}.
\]
Notice that for this collection of elements in $\Phi$, we do not require they represent points in the sense of Definition~\ref{points}, although in some cases (such as generally prime), it is superfluous.

Let $K_p^{\wedge}(C)=\{C'\in K_p^{\wedge}(\Phi)~|~C'\geq C\}$ for $C\in \Phi$. Again we need a formal empty set to be included in this collection to give a topology. For a homomorphism $f:\Phi\rightarrow \Psi$ of complete lattices, define
$K_p^{\wedge}(f):K_p^{\wedge}(\Psi)\rightarrow K_p^{\wedge}(\Phi)$ as a supremum
\[
K_p^{\wedge}(f)(C)=\bigvee f^{-1}(\downarrow(C)),
\]
where $\downarrow(C)$ denotes the \textit{principal ideal} generated by $C$ meaning that it contains all elements below $C$, for every meet prime $C$ in $\Psi$. Such $K_p^{\wedge}(f)$ is well-defined. Indeed, suppose $\bigvee f^{-1}(\downarrow(C))\geq A\wedge B$. Then by applying $f$, we have $C= \bigvee (\downarrow(C))\geq f(A)\wedge f(B)$, which implies $C\geq f(A)$ or $C\geq f(B)$ since $C$ is meet prime. We assume $C\geq f(A)$ without loss of generality. It follows that $\bigvee f^{-1}(\downarrow(C))=\bigvee_{D\in f^{-1}(\downarrow(C))}D=\bigvee_{C\geq f(D)} D\geq A$.

\begin{prop}\label{ptspace}
For a complete distributive lattice $\Phi$, the collection $\{K_p^{\wedge}(C)~|~C\in\Phi\}$ forms a topology of closed subsets of $K_p^{\wedge}(\Phi)$, and $K_p^{\wedge}$ is a contravariant functor.
\end{prop}

\pf This is essentially Lemma II.1.3 and Theorem II.1.4 in~\cite{JohnstoneP}. For topology, we have a formal empty set, and $K_p^{\wedge}(0)=K_p^{\wedge}(\Phi)$. It is also clear that $K_p^{\wedge}(\bigwedge C_i)=\bigwedge K_p^{\wedge}(C_i)$. We need to show $K_p^{\wedge}(C_1)\cup K_p^{\wedge}(C_2)=K_p^{\wedge}(C_1\wedge C_2)$ for each $C_i\in\Phi$. Notice that $K_p^{\wedge}(C_i)\subseteq K_p^{\wedge}(C_1\wedge C_2)$ by definition, so that $K_p^{\wedge}(C_1)\cup K_p^{\wedge}(C_2)\subseteq K_p^{\wedge}(C_1\vee C_2)$. Conversely, if $C\geq C_1\wedge C_2$, then $C\geq C_i$ for some $i$ by primeness, which implies $K_p^{\wedge}(C_1)\cup K_p^{\wedge}(C_2)\supseteq K_p^{\wedge}(C_1\vee C_2)$. Furthermore, for any homomorphism $f:\Phi\rightarrow\Psi$, we claim that $K_p^{\wedge}(f)$ is continuous. Indeed,
\[
K_p^{\wedge}(f)^{-1}(K_p^{\wedge}(C))=K_p^{\wedge}(f(C)),
\]
which follows immediately from the fact that $B\geq f(A)$ if and only if $K_p^{\wedge}(f)^{-1}(B)\geq A$, where $A,B\in\Phi$. Finally, for an identity map $id$ of lattices, it is clear that $K^{\wedge}_p(id)(C)=C$ for any $C$. Next suppose $f,g$ are lattice maps and $g\circ f$ is their composition, by exhibiting the index sets we have \[
K^{\wedge}_p(f)\circ K^{\wedge}_p(g)(C)=\bigvee_{\bigvee\limits_{C\geq g(D)}D ~\geq~ f(E) }E=\bigvee_{C\geq gf(E)}E=K^{\wedge}_p(g\circ f)(C),
\]
as required.~~~\eop



Following the same pattern, we define spaces consisting of completely meet prime elements, (completely) meet irreducible elements, and their dual cases. Denoted by $K^{\wedge}_{cp}(\Phi)$, $K^{\wedge}_i(\Phi)$, $K^{\wedge}_{ci}(\Phi)$ and dually $K^{\vee}_p(\Phi)$, $K^{\vee}_{cp}(\Phi)$, $K^{\vee}_{i}(\Phi)$ and $K^{\vee}_{ci}(\Phi)$, respectively. With a slightly different argument, we can show similar results to Proposition~\ref{ptspace}.

Notice that in the case of completely join primes, this space coincides with our $K_{gp}(\Phi)=K^{\vee}_{cp}(\Phi)$ since every such prime element defines a point in $K(\Phi)$ (see the argument above Theorem~\ref{closedprime}). Also, under the dualization of a lattice $\Phi$, $K_p(\Phi)$ becomes a subspace of $K_p^{\vee}(\Phi)$ which consists of all join prime elements in $\Phi$. However, in $K_p(\Phi)$ we require additionally that each prime $C$ also defines a point in $K(\Phi)$.

Similar to Theorem~\ref{closedprime} (which is in fact the case of $K^{\vee}_{cp}(\Phi)$) and its proof, we have a family of parallel results given in Theorem~\ref{classsification} and its corollary.

\begin{thm}\label{classsification}
Let $\Phi$ be a complete lattice of subcategories of certain type in a category $\al$. Then the subcategories generated by join primes are classified by the closed subsets of $K^{\vee}_p(\Phi)$. More precisely, there is a bijection
\[
\th:\{K^{\vee}_p(C)\}_{C\in\Phi}\rightarrow\{C\in\Phi~|~C=\widehat{C}\}:\xi
\]
where $\widehat{C}=\bigvee_{\substack{C'\leq C\\ C'~\text{join prime}}}C'$, defined by $\th(K^{\vee}_p(C))=\widehat{C}$ and $\xi(C)=K^{\vee}_p(C)$.
\end{thm}

\pf It is straightforward to check that $\th$ and $\xi$ are well-defined and inverses to each other, noting that $\widehat{\widehat{C}}=\widehat{C}$ and $K^{\vee}_p(\widehat{C})=K^{\vee}_p(C)$.~~~\eop

With a slight modification of the proof, we can deduce similar results to Theorem~\ref{classsification}.

\begin{coro}\label{cortoclassification}
Let $\Phi$ be a complete lattice of subcategories of certain type in a category $\al$. Then the subcategories generated by (complete resp.) join irreducibles are classified by the closed subsets of $K^{\vee}_i(\Phi)$ ($K^{\vee}_{ci}(\Phi)$ resp.). The dual statements for the meet primes and irreducibles also hold accordingly.
\end{coro}


\section{Comparing classical results}\label{cmpclassical}


There are several examples that fit into the context of Theorem~\ref{closedprime} in the sense that each subcategory is g-primely generated. The cases when we already know the generators would be easier to handle; for example, Serre subcategories in a noetherian abelian category, thick subcategories in $D_{perf}(R)$ and localizing subcategories in a derived category $D(R)$ of a commutative noetherian ring $R$, or more generally, localizing subcategories in a stable homotopy category.

The proofs are in fact restated from the known results, but they are reorganized for our purpose to show that the subcategories are g-primely generated, so that Theorem~\ref{closedprime} applies to give a classification by our classifying spaces. To avoid repetition we only include the last example as a demonstration. As usual, we use the brackets $\langle-\rangle$ to denote a subcategory of certain type generated by a collection of objects.

\begin{itemize}

\item \textit{Localizing subcategories in a stable homotopy category}

\end{itemize}

A \textit{stable homotopy category} $\al$ is a triangulated category with a compatible closed symmetric monoidal structure, in which arbitrary coproducts exist, every cohomology functor is representable and there is a set of strongly dualizable objects generating $\al$. See \cite{HoveyPS} for detail. The following lemma gives a description of localizing subcategories in $D(R)$, which has a parallel result in a more general setting of noetherian stable homotopy category.

\begin{lem}
Let $T\subseteq \Spec R$ be an arbitrary subset. Then we can characterize the localizing subcategory generated by the residue field $k(\mp)$ parametrized by $T$ as
\[
\langle k(\mp)~|~\mp\in T\rangle=\{X\in D(R)~|~X\otimes K_T=0 \}
\]
where $K_T=\bigoplus_{\mp\in T^c}k(\mp)$ is a coproduct over the compliment $T^c$.
\end{lem}

\pf Since $k(\mp)\otimes k(\mq)=0$ if $\mq\neq\mp$ and that tensor product respects coproducts, one containment holds. Now suppose $X\otimes K_T=0$ so that $X\otimes k(\mp)=0$ for every $\mp\in T^c$. Thus by tensoring the complex of injectives supported at a closed point,  $\Gamma_{\overline{\mp}/\overline{\mp}-\{\mp\}}(X)\otimes k(\mp) =X\otimes k(\mp)=0$, hence $\Gamma_{\overline{\mp}/\overline{\mp}-\{\mp\}}(X)=0$ by Lemma 2.14 in~\cite{Neeman}. Therefore, $X\in\langle k(\mp)~|~\mp\in T\rangle$ by Lemma 2.10 in~\cite{Neeman}.~~~\eop


Recall from Section 6 of~\cite{HoveyPS} that a \textit{noetherian stable homotopy category} is a monogenic stable homotopy category (i.e. the generating set consists of only the unit object which is also compact) with a single compact generator $S$ such that the graded commutative ring $\pi_\ast(S)=[S,S]_\ast$ is noetherian, denoted by $R$. The derived category of a commutative noetherian ring gives such an example. For any prime $\mp$ of $R$, denote $K(\mp)=S_\mp\wedge S/\mp$, where $S_\mp$ is the localization at $\mp$ and $S/\mp$ is the Koszul complex associated to $\mp$.

\begin{lem}\label{lcznshc}
Let $\al$ be a noetherian stable homotopy category such that every monochromatic category $\langle K(\mp)\rangle$ is minimal among nonzero localizing subcategories of $\al$. Then every localizing subcategory has the form
\[
\langle K(\mp)~|~\mp\in T\rangle=\{X\in\al~|~X\wedge K_T=0\}
\]
for some subset $T\subseteq \Spec R$, where $K_T=\coprod_{\mp\in T^c}K(\mp)$.
\end{lem}

\pf This is in fact Corollary 6.3.4 in~\cite{HoveyPS}. Since $K(\mp)\wedge K_T=0$ for any $\mp\in T^c$ by Proposition 6.1.7 in~\cite{HoveyPS}, we have $\langle K(\mp)~|~\mp\in T\rangle\subseteq\{X\in\al~|~X\wedge K_T=0\}$. Conversely, suppose $X\in\al$ such that $X\wedge K_T=0$. Notice that $X\wedge K_T=0$ if and only if $X\wedge K(\mp)=0$ for every $\mp\in T^c$. It follows that $X\in\langle X\wedge K(\mp)~|~\mp\in\Spec R\rangle$ by Proposition 6.3.2 in~\cite{HoveyPS}, so that $X\in\langle X \wedge K(\mp)~|~\mp\in T\rangle\subseteq\langle K(\mp)~|~\mp\in T\rangle$ since every localizing subcategory is an ideal.~~~\eop

\begin{prop}\label{Nstabhcat}
Let $\al$ be a noetherian stable homotopy category such that every monochromatic category $\langle K(\mp)\rangle$ is minimal among nonzero localizing subcategories of $\al$. Then every localizing subcategory is g-primely generated.
\end{prop}

\pf Thanks to Lemma~\ref{lcznshc}, we assume that $\langle K(\mq)~|~\mq\in T_i\rangle$ for $i\in I$ are localizing subcategories. Then $K(\mp)\in \bigvee_{i\in I}\langle K(\mq)~|~\mq\in T_i\rangle$ implies
\begin{eqnarray*}
K(\mp)&\in &\bigvee_{i\in I}\langle K(\mq)~|~\mq\in T_i\rangle=\langle K(\mq)~|~\mq\in \bigcup_{i\in I}T_i\rangle\\
&=&\{X\in\al~|~X\wedge K_{\bigcup_{i\in I}T_i}=0\}.
\end{eqnarray*}
Hence $K(\mp)\wedge K_{\bigcup_{i\in I}T_i}=0$, so that $\mp\in \bigcup_{i\in I}T_i$ and $\mp\in T_i$ for some $i\in I$. In particular, $K(\mp)\in\langle K(\mq)~|~\mq\in T_i\rangle$ for some $i$.~~~\eop

\section{An example}\label{exceptional}
Nevertheless, we do have examples of subcategories that are not g-primely generated. Consider the quiver $A_2:\stackrel{1}{\circ}\longleftarrow\stackrel{2}{\circ}$ and the category $\al$ of finite dimensional representations of $A_2$ over a field $k$. For convenience, we denote $a=P_1$, $b=P_2$ and $c=S_2$.


Let $\Phi$ be the lattice of nullity classes (subcategories that are closed under quotients and extensions) in $\al$. We can specify $\Phi$, $K(\Phi)$, $K_p(\Phi)$, $K_{gp}(\Phi)$ and their topology $K_p(C)$ explicitly thanks to the Auslander-Reiten quiver of $A_2$
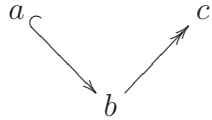
\[
\xymatrix{
a\ar@{^{(}->}[dr]&&c\\
&b\ar@{->>}[ur]&
}
\]
Thus
$\Phi=\{\langle 0\rangle, \langle a\rangle,\langle c\rangle,\langle b,c\rangle,\al\}$ and
$K(\Phi)=\{P_{\langle 0\rangle},P_{\langle a\rangle},P_{\langle c\rangle},P_{\langle b,c\rangle}\}$, so that $K_{gp}(\Phi)=K_p(\Phi)=\{P_{\langle 0\rangle},P_{\langle a\rangle},P_{\langle c\rangle}\}$, with 5 closed subsets $\emptyset$, $\{P_{\langle0\rangle}\}$, $\{P_{\langle0\rangle},P_{\langle a\rangle}\}$, $\{P_{\langle0\rangle},P_{\langle c\rangle}\}$ and $K_{gp}(\Phi)$. This example demonstrates that nullity class cannot be classified by a topological space but only those primely generated ones. Notice that $P_{\langle b,c\rangle}$ is not prime since $\langle b,c\rangle\leq\langle a\rangle\vee\langle c\rangle$.

However, nullity class is classified by closed subsets with extra conditions in a different topological space. For example, in the spectrum $\Spec\al=\{a,b,c\}$ of premonoforms, they correspond to the closed and extension-closed subsets, see Theorem 6.8 in~\cite{liu2}.

Finally, it is interesting to point out that the lattice $\Phi$ does not come from a topological space by computing its \textit{point space} $\mathrm{pt}\Phi$ (see the definition of point space in Chapter II.1 of~\cite{JohnstoneP}), on which there is no valid topology. As a matter of fact, by depicting the diagram of the lattice $\Phi$
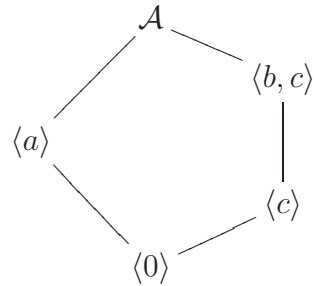
\[
\xymatrix@R=0.5em{
&\al\ar@{-}[dr]&\\
&&\langle b,c\rangle\ar@{-}[dd]\\
\langle a\rangle\ar@{-}[uur]\ar@{-}[ddr]&&\\
&&\langle c\rangle\\
&\langle0\rangle\ar@{-}[ur]&
}
\]
and a general theory of lattices (see e.g. Chapter 4 of~\cite{DaveyPriestley}), we conclude that the lattice $\Phi$ is not even distributive since it contains a pentagon, see Theorem 4.10 (ii) in~\cite{DaveyPriestley}.

\bibliography{mybib}

\begin{thebibliography}{10}

\bibitem{Balmer}
P.~Balmer.
\newblock The spectrum of prime ideals in tensor triangulated categories.
\newblock {\em Journal f\"{u}r die reine und angewandte Mathematik},
  (588):149--168, 2005.

\bibitem{BBD}
A.~A. Beilinson, J.~Bernstein, and P.~Deligne.
\newblock Faiscearx pervers.
\newblock {\em Asterisque}, (100):5--171, 1982.

\bibitem{Benson}
D.~Benson, S.~B. Iyengar, and H.~Krause.
\newblock Statifying modular representations of finite groups.
\newblock {\em Annals of Mathematics}, (174):1643--1684, 2011.

\bibitem{DaveyPriestley}
B.~A. Davey and H.~A. Priestley.
\newblock {\em Introduction to Lattices and Order}.
\newblock Cambridge University Press, 2002.

\bibitem{Gabriel}
P.~Gabriel.
\newblock Des cat\'{e}gories ab\'{e}liennes.
\newblock {\em Bulletin de la Soci\'{e}t\'{e} Math\'{e}matique de France},
  (90):323--448, 1962.

\bibitem{Hochster}
M.~Hochster.
\newblock Prime ideal structure in commutative rings.
\newblock {\em Transactions of the American Mathematical Society},
  (142):43--60, 1969.

\bibitem{HoveyPS}
M.~Hovey, J.~H. Palmieri, and N.~P. Strickland.
\newblock {\em Axiomatic Stable Homotopy Theory}, volume 128 of {\em Memoirs of
  the American Mathematical Society}.
\newblock American Mathematical Society, 1997.

\bibitem{JohnstoneP}
P.~Johnstone.
\newblock {\em Stone Spaces}, volume~3 of {\em Cambridge Studies in Advanced
  Mathematics}.
\newblock Cambridge University Press, 1982.

\bibitem{Kanda}
R.~Kanda.
\newblock Classifying serre subcategories via atom spectrum.
\newblock {\em Advances in Mathematics}, (231):1572--1588, 2012.

\bibitem{Kock}
J.~Kock and W.~Pitsch.
\newblock Hochster duality in derived categories and point-free reconstruction
  of schemes.
\newblock arXiv:1305.1503, May 2013.

\bibitem{liu2}
Y.~Liu and D.~Stanley.
\newblock A classification of torsion classes in abelian categories.
\newblock arXiv:1610.09528, Feburary 2017.

\bibitem{Neeman}
A.~Neeman.
\newblock The chromatic tower for {$D(R)$}.
\newblock {\em Topology}, (3):519--532, 1992.

\bibitem{Neeman2}
A.~Neeman.
\newblock Colocalizing subcategorie of {$D(R)$}.
\newblock {\em Journal f\"{u}r die reine und angewandte Mathematik},
  (653):221--243, 2011.

\bibitem{Rosenberg}
A.~Rosenberg.
\newblock Reconstruction of schemes.
\newblock Max Planck Institute, 1996.

\bibitem{Rouquier}
R.~Rouquier.
\newblock Derived categories and algebraic geometry.
\newblock In {\em Triangulated Categories}, volume 375 of {\em London
  Mathematical Society Lecture Note Series}. Cambridge University Press, 2010.

\bibitem{Stan}
D.~Stanley.
\newblock Invariants of $t$-structures and classification of nullity classes.
\newblock {\em Advances in Mathematics}, (224):2662--2689, 2010.

\bibitem{Alonso2}
L.~Alonso Tarrio, A.~Jeremias Lopez, and M.~Saorin.
\newblock Compactly generated $t$-structures on the derived category of a
  noetherian ring.
\newblock {\em Journal of Algebra}, (324):313--346, 2010.

\end{thebibliography}
\bibliographystyle{plain}

\end{document}